\documentstyle[amstex,amscd,12pt]{article}
\makeatletter
\date{}
\newtheorem{theorem}{Theorem}[section]

\newtheorem{proposition}[theorem]{Proposition}
\newtheorem{definition}[theorem]{Definition}

\newtheorem{problem}[theorem]{Problem}
\newcommand{\edim}{{\rm e-dim}}
\newcommand{\z}{{\Bbb Z}}

\newcommand{\lo}{\longrightarrow}
\newcommand{\sm}{\setminus}

\begin{document}

\title{Some mapping theorems for extensional dimension\\
}

\author{Michael  Levin  and Wayne Lewis}

\maketitle
\begin{abstract} We present some results related to theorems 
of Pasynkov and Torunczyk on the geometry of maps of
finite dimensional compacta. 
\bigskip
\\
{\bf Math. Subj. Class.:} 55M10, 54F45.
\end{abstract}
\begin{section}{Introduction} 

All topological spaces are assumed to be separable metrizable, $I=[0,1]$.
Recall that the covering dimension  $\dim X$ is the smallest number $n$
such that every open cover of $X$ admits an open refinement of order 
$\leq n+1$.  The cohomological dimension $\dim_G X$ with respect to an abelian
group  $G$ is  the smallest number $n$ such that ${\check H}^{n+1}(X,A;G)=0$
for all closed subsets $A$ of $X$. 

By  a classical result of Alexandroff 
$\dim_\z  X =\dim X$  for a finite dimensional $X$. Solving an outstanding
problem in Dimension Theory  Dranishnikov constructed in 1986 
an infinite dimensional compactum (=compact metric space) of $\dim_\z =3$.
Many classical results for the covering  and the cohomological dimensions
shows 
a great deal of similarity between the theories.
Here we  mention just two such results.

\begin{theorem}(Hurewicz \cite{En})
\label{i1}
Let $ f: X \lo Y$ be a map of compacta. Then $\dim X \leq 
\dim Y +\dim  f$, where $\dim f = \sup \{ \dim f^{-1}(y): y \in Y \}$.
\end{theorem}

\begin{theorem}(Skljarenko \cite{br})
\label{i2}
Let $ f: X \lo Y$ be a map of compacta. Then $\dim_\z X \leq 
\dim_\z Y +\dim_\z  f$, 
where $\dim_\z f = \sup \{ \dim_\z f^{-1}(y): y \in Y \}$.
\end{theorem}

Despite such similarities
 the covering dimension was   investigated
using  mostly set-theoretic methods and the cohomological
dimension  was investigated using mostly  algebraic methods.
In general the results for  the cohomological dimension were proved
with difficulty exceeding   similar results for the covering dimension. 
 The 
proofs of the theorems 
mentioned above   are  a  good illustration of the differences in  
the approaches.
 The essentially  different definitions of the dimensions
seem to  justify    using different tools. 
However the following  two  theorems offer a common point of view on 
the theories.

 \begin{theorem}(Hurewicz-Wallman \cite{En}) \label{i3}
$\dim X\le n$ if and only if 
 any map $f:A\to S^{n}$ of a closed subset $A$ of $X$  into the $n$-dimensional
 sphere $S^n$ 
extends over $X$.
\end{theorem}

\begin{theorem}(Cohen {\cite{Co}})
\label{i4}
$\dim _{G } X\le n$ 
if and only if 
 any map $f:A\to K(G,n)$ of a closed subset $A$ of $X$  into 
 an Eilenberg-MacLane complex   of type $(G,n)$
extends over $X$.
\end{theorem}

Theorems \ref{i3} and \ref{i4}
 suggest that one may try to consider extension properties of 
maps into different CW-complexes hoping to create a theory 
encompassing both the covering and the cohomological dimension theories.
Indeed,  in
the early 1990's  Dranishnikov laid down a foundation of such theory which
he called ``Extension  Theory".  Extension Theory not only provided 
simpler proofs for many results in cohomological dimension theory, but
it also gave an explanation of certain phenomena in cohomological dimension
which  are sometimes more general than expected, see \cite{drdyd} and
\cite{drdyd1}.

Let us introduce the terminology.
The extensional dimension of  $X$ 
does not
exceed a CW-complex $K$, written e-dim$X \leq K$, if any map of a closed
subset of $X$ into $K$ extends over $X$. Thus for the covering dimension
dim$X \leq n$ if and only if e-dim$X \leq S^n$
and for the cohomological dimension 
dim$_G X \leq n$ if and only if e-dim$X \leq K(G,n)$.
 Let $f : X \lo Y$ be a map.
We will use the following notation:
 
  e-dim$f \leq K$
if e-dim$f^{-1}(y) \leq K$ for every $y \in Y$;

e-dim$Z \leq $e-dim$f$
if for every CW-complex $K$,
e-dim$f \leq K$ implies e-dim$Z \leq K$;
 
e-dim$X \leq $e-dim$Y$  if 
 for every CW-complex $K$, e-dim$Y \leq  K$  implies e-dim$X \leq K$;
 
 for CW-complexes $K$ and $L$, $K \leq L$ if for every  $X$,
 e-dim$X \leq K$ implies e-dim$X \leq L$. 
 
 This note is mainly devoted to  an extensional dimension
  generalization  of some results
 of Pasynkov \cite{p} and Torunczyk \cite{t} on the geometry of maps.
 For example,  we  obtain the following  
 versions of   Theorems \ref{i1} and \ref{i2}. Note that Theorem \ref{2t1}
 is an improvement  of the results of \cite{l3} which provides
 the  right  estimates 
for cohomological dimension.

   \begin{theorem}
\label{2t1}
(cf. Remark 2, Section 5)
Let $f : X \lo Y$ be a map of compacta and let CW-complexes  $K$ and $L$
be such that $K$ is countable, $ \edim f\leq K$ and
$\edim Y \leq L$. Then $X_* =X \times I$ can be decomposed into
$X_* = A\cup B$ for subspaces $A$ and $B$ such that 
$\edim A \leq K$ and $\edim B \leq L$. As a result  of
such a decomposition  $\edim X_* \leq K*L$.
\end{theorem}

\begin{theorem}
\label{smash}
Let $f : X \lo Y$ be a map of compacta with
$X$ finite dimensional and let $K$ and $L$ be
CW-complexes such that $K$ is countable,
 $\edim f \leq K$ and $\edim Y \leq L$.
Then $\edim X \leq  K \wedge L$.
\end{theorem}

Pasynkov \cite{p} and Torunczyk \cite{t} proved the following
remarkable theorem.
\begin{theorem}
\label{pt}
(Pasynkov-Torunczyk) Let $f : X \lo Y$ be a map
of finite dimensional compacta $X$ and $Y$.
Then the following conditions are equivalent:

(1) $\dim f \leq n$;

(2) there exists a $\sigma$-compact subset $A$ of $X$ such that
$\dim A \leq n-1$ and $\dim f|_{X\sm A} \leq 0$;

(3) almost every  map $g: X \lo I^n$ has the property that the
map $(f,g) : X \lo Y \times I^n$ is $0$-dimensional, where
almost=all but a set of first category;

(3') there exists a map $g: X \lo I^n$ such that the map
$(f,g) : X \lo Y \times I^n$ is $0$-dimensional.
\end{theorem}

In this note we will prove two theorems generalizing Theorem \ref{pt}.
\begin{theorem}
\label{t1}
Let $f : X \lo Y$ be a map of compacta and let $K$ be a CW-complex.
Consider the following properties:

(1) $ \edim f \leq \Sigma K$;

(2) there exists a $\sigma$-compact subset $A$ of $X$ such that $\edim A\leq K$
and $\dim f|_{X\sm A} \leq 0$;

(3)  almost every map $g : X \lo I$ is such that for the map $(f,g) : X \lo Y\times I$,
$\edim (f,g) \leq K$;

(3') there exists a map $g : X \lo I$ such that  for the map
$(f,g) : X \lo Y\times I, \edim (f,g) \leq K$.

Then:
(3)$\Rightarrow$(2)$\Rightarrow$(1),
 (3)$\Rightarrow$(3')$\Rightarrow$(1) and if 
  $Y$ is finite-dimensional and $K$ is countable then
  (1)$\Rightarrow$(3).

In particular all the properties are equivalent if $Y$ is finite dimensional 
and $K$ is countable.
\end{theorem}

Note that 
Theorem \ref{pt} can be derived from Theorem \ref{t1}. Indeed,
it is obvious for the properties (1) and (2) of Theorem \ref{pt}.
One can easily show by induction that the property (3) of Theorem \ref{t1}
implies that the set of functions $g$ satisfying the property  (3) of Theorem \ref{pt}
for $\dim f \leq n$ is dense in $C(X,I^n)$. On the other hand this set is
$G_\delta$ in $C(X,I^n)$ and hence  (3) holds for $\dim f \leq n$. Clearly
(3)$\Rightarrow$(3')$\Rightarrow$(1) and Theorem \ref{pt} follows
(see also Remark 1, Section 5).

It is unknown if the finite dimensional restriction on $Y$ is necessary in
Theorems \ref{pt} and \ref{t1}. Some versions of these theorems without this 
restriction were obtained in \cite{st,l1,l3}.  The following theorem improves
on these results.
\begin{theorem}
\label{t2}
Let $f : X \lo Y$ be a map of compacta and let $K$ be a CW-complex.
Define  $f_* : X_*=X \times I \lo Y$ by $f_*(x,t)=f(x)$.
Consider the following properties:

(0) $\edim f\leq K$;

(1) $\edim f_* \leq \Sigma K$;

(2) there exists a $\sigma$-compact subset $A$ of $X_*$ such that
 $ \edim A\leq K$
and $\dim f_* |_{X_{*}\sm A} \leq 0$;

(2') a set $A$  satisfying (2) can be chosen such that $A$ splits
into compacta $A=\cup_{i=1,\infty} A_i$  such that
each component of $A_i$ admits a $0$-dimensional
map into a fiber of $f$.

(3)  almost every map $g_* : X_* \lo I$
is such that for the map $(f_*,g_*) : X_*\lo Y\times I$,
$\edim (f_*,g_*) \leq K$;

Then (0)$\Rightarrow$(3)$\Rightarrow$(2')$\Rightarrow$(2)$\Rightarrow$(1).

\end{theorem}

Theorems \ref{t1} and \ref{t2} are proved in Section 4.  Applications
of Theorem \ref{t2} including the proofs of Theorems \ref{2t1} and
 \ref{smash} 
 are presented in Section 2. Section 3 is devoted to 
Krasinkiewicz  maps
which are essentially used  in the proofs. In Section 5 we discuss
some  related results and problems.

Finally we wish to thank the referee  for pointing  out that Proposition \ref{pp3} can also be
derived  from Lemma 2  of    \cite{t}

\end{section}

\begin{section}{ Applications}
We will use the following facts.
\begin{theorem}
\label{dru} (\cite{dru})
Let $f : X \lo Y$ be a light map of compacta.  Then $\edim X \leq  \edim Y$.
\end{theorem}
In some cases we will need a  stronger result than Theorem \ref{dru}.
\begin{definition}
(\cite{kat})
A map $f : X  \lo Y$ of metric spaces $X$ and $Y$ is said to be
uniformly $0$-dimensional if there is a metric on $X$
such that for every $\epsilon>0$
 every point in $f(X)$ 
has a neighborhood $V$ in $Y$ such that $f^{-1}(V)$ splits into
a disjoint family of open sets of diam$ < \epsilon$.
\end{definition}
Examples of uniformly $0$-dimensional maps are: a $0$-dimensional map
of compacta, a $0$-dimensional perfect (=closed with compact fibers) map
of metric spaces, a map of compacta restricted to the union of all
trivial components of its fibers.  

\begin{theorem}
\label{l5}
(\cite{l5})
\label{unif}
Let $f :X \lo Y$ be a uniformly $0$-dimensional map of metric spaces.
Then $\edim X \leq \edim Y$.
\end{theorem}

\begin{theorem}
\label{ol}
(\cite{ol})
Let $K$ be a countable CW-complex and let $A$ be a subspace of a compactum $X$
such $\edim A \leq K$. Then there is a completion $A' \subset X$ of $A$ 
such that $\edim A' \leq K$.
\end{theorem}

\begin{theorem}
\label{dyd}
(\cite{dyd}, see also Remark 2, Section 5) Let $K$ and $L$ be CW-complexes and let $X=A \cup B$ be 
a decomposition of a separable metric space $X$ into subspaces $A$ and $B$
such that $\edim A\leq K$ and  $\edim B \leq L$. Then $\edim X \leq K*L$.
\end{theorem}

\begin{theorem}
\label{rub}
(\cite{rub},see also \cite{dyd})
Let $X=A\cup B$ be a decomposition of a separable  metric space $X$.
Then $\dim_\z X \leq \dim_\z A +\dim_\z B +1$.
\end{theorem}

\begin{theorem} 
\label{dr1}
(\cite{drdyd})
Let $X$ be a  compactum and let
$K$ be a simply connected CW-complex. Consider the following
conditions:

(1) $\edim X \leq  K$;

(2) $\dim_{H_i(K)} X \leq i$ for every $i > 1$;
 
(3) $\dim_{\pi_i(K)} X \leq i$ for every $i >1$.

Then   (2) and (3) are equivalent and (1) implies both (2) and (3).
If $X$ is  finite dimensional then all the conditions are equivalent. 
 In particular if $X$ is  finite dimensional then
 $\dim_G X \leq n$, $n >1$ if and only if
 $\edim X \leq M(G,n)$ where $M(G,n)$ is a Moore space of type $(G,n)$.
\end{theorem}

\begin{theorem}
\label{dr2}
(\cite{dr}) Assume that a compactum $X$ is expressed as the union
$X=A \cup B$. Then $\dim_\z X \leq \dim_\z (A \times B)$ if 
$\dim_\z X^2 =2\dim_\z X$
and $\dim_\z X \leq \dim_\z (A \times B)+1$  if
$\dim_\z X^2 \neq 2\dim_\z X$ (that is $\dim_\z X^2 =2\dim_\z X-1$).
\end{theorem}

\begin{proposition}
\label{2p0}
(\cite{l3})
Let $K$ be a CW-complex and let each component of a compactum $X$
be of $\edim \leq K$. Then $\edim  X \leq K$.
\end{proposition}

The following two proposition can be easily derived from
the proofs of Proposition 12, \cite{dr}
 and Lemma 2.11, \cite{drrepsh} respectively.

\begin{proposition}
\label{2p1}
Let $X$ and $Y$
 be compacta and let $A \subset X$. Then
there is a completion $A' \subset X$ of $A$
such that $\dim_\z (A' \times Y)=\dim_\z
(A \times Y)$.
\end{proposition}

\begin{proposition}
\label{2p2}
Let $K$ be a CW-complex such that $\edim (X\times Y) \leq K$ where $X$ and $Y$
are non-trivial continua. Then $K$ is simply connected.
\end{proposition}

It was shown in \cite{l3} that for a map of compacta $f : X\lo Y$ and
CW-complexes $K$ and $L$,
$X$ can be decomposed into $X =A \cup B$ such that e-dim$A \leq K$ and
$\edim B \leq L$ 
provided e-dim$f \leq K$, e-dim$Y \leq L$ and $K$ is countable. 
Theorem \ref{2t1}  improves this result.

{\bf Proof of Theorem \ref{2t1}.}
 By Theorem  \ref{t2}
there is $A \subset X_*$ such that $\edim A \leq K$ and
$\dim f_* |_B \leq 0$ where $B = X_* \sm A$.
By  Theorem \ref{ol} there is
a completion $A' \subset X_*$ of $A$ with e-dim$A' \leq K$.
$B' =X_* \sm A'  \subset B$ is $\sigma$-compact and hence by Theorem \ref{dru},
e-dim$B' \leq L$. By Theorem \ref{dyd}
$\edim X \leq K*L$.
\hfill $\Box$
\\

 Note that applying  Theorem \ref{rub} to 
  the decomposition of $X_*$ from Theorem \ref{2t1} for
  $\dim_\z$ it follows that
$\dim_\z X +1 =\dim_\z X_* =
\dim_\z  (A' \cup B') \leq  \dim_\z A' +\dim_\z  B' +1$ and
we have obtained  Theorem  \ref{i2},
see \cite{br} for related results.

Dranishnikov and Dydak \cite{drdyd1}  proved that if for 
finite dimensional compacta $X$ and $Y$
and CW-complexes $K$ and $L$, $\edim X \leq K$ and $\edim Y \leq L$
then $\edim (X \times Y) \leq K \wedge L$. Theorem \ref{smash}
is a mapping version of Dranishnikov-Dydak's result.

{\bf Proof of Theorem \ref{smash}.} If $K$   is disconnected then $\edim f =0$  and 
if $L$ is disconnected then $\dim Y =0$. These cases can be treated
in a way similar to the one in the proof of Theorem 5.6 of \cite{drdyd1}
to show that $\edim X \leq K \wedge L$. 
 
Assume that both $K$ and $L$
are connected. Then $K\wedge L$ is simply connected.
By Theorem \ref{2t1}, $\edim (X_*) \leq K*L$, $X_* =X\times I$
 and since  $K*L$ is homotopy equivalent
to $\Sigma (K \wedge L)$,    $\edim (X_*) \leq    \Sigma (K \wedge L)$.
 Then
 by Theorem
\ref{dr1}, $\dim_{H_{i+1}(\Sigma (K \wedge L))} (X_*)=
\dim_{H_{i} ( K \wedge L)} (X_*)=
\dim_{H_{i}(K \wedge L)} (X) +1  \leq i+1$.
Hence $\dim_{H_{i}(K \wedge L)} (X)  \leq i$ and again by 
Theorem \ref{dr1}, $\edim X \leq K \wedge L$.
\hfill $\Box$\\\\
Note that if in the proof of Theorem \ref{smash}
 one uses   Theorem \ref{5t2}  (see
Section 5, Remark  2 of this note) instead of 
Theorem \ref{2t1} then the restriction in  
Theorem \ref{smash} that $K$ is countable
can be omitted.

The following   theorem  was 
 announced  by  Dranishnikov in 1996. We will present a proof of this
 result based on Theorem \ref{t2}.

\begin{theorem}
\label{t3}
(Dranishnikov) 
Let $f: X \lo Y$ be a map of compacta. Then:
$\dim_\z X \leq$${\rm sup}\{ \dim_\z$ $(f^{-1}(y)\times Y):y\in Y \}$
if $\dim_\z X^2 =2\dim_\z X$
and $\dim_\z X \leq {\rm sup}\{ \dim_\z (f^{-1}(y)\times Y):y\in Y \}+1$
otherwise.
\end{theorem}
{\bf Proof.}
Let  $n=$sup$\{ \dim_\z (f^{-1}(y)\times Y):y\in Y \}$.
Take a $\sigma$-compact $A \subset X_*$ satisfying (2') of Theorem \ref{t2},
that is $A$ splits into compacta $A=\cup A_i$ such that
each component of each $A_j$
admits a $0$-dimensional map to a fiber of $f$. Then each component
of $A_i \times Y$ admits a $0$-dimensional map to $f^{-1}(y) \times Y$
for some $y \in Y$ and hence 
by Theorem \ref{dru} and Proposition \ref{2p0}
$\dim_\z (A \times Y) \leq n$ and 
by Proposition \ref{2p1} there is a completion $A'\subset X_*$
of $A$ such that $\dim_\z (A' \times Y) \leq n$. Then $B'=X_* \sm A'$
is $\sigma$-compact and hence splits into 
compacta $B'=\cup B_j$ admitting $0$-dimensional maps
into $Y$. It follows that $A' \times B_j$ admits
a $0$-dimensional perfect map to $A'\times Y$. Thus
 by Theorem \ref{unif} $\dim_\z (A' \times B_j) \leq n$
 and hence 
$\dim_\z(A' \times B') \leq n$. By Theorem \ref{dr2}
$\dim_\z X_* \leq n +1$ if $\dim_\z X_*^2=2\dim_\z X^*$ and 
$\dim_\z X_* \leq n +2$  otherwise.
Clearly $\dim_\z X_*=\dim_\z X +1$ and $\dim_\z X_*^2=\dim_\z X^2 +2$ and the theorem
follows.
\hfill $\Box$
\\

Note that it is unknown if 
the  case $\dim_\z X =\dim_\z A + \dim_\z B +2$ 
in Theorem \ref{dr2} ever occurs.
Eliminating this case would result in eliminating 
the case
$\dim_\z X  = {\rm sup}\{ \dim_\z (f^{-1}(y)\times Y):y\in Y \}+1$
in Theorem \ref{t3}.

Our next application is to  prove in a slightly
stronger form  the generalized Hurewicz formula  obtained by
Dranishnikov, Repov\v s and {\v S}{\v c}epin \cite{drrepsh}.

\begin{theorem}
\label{t4}
(Dranishnikov-Repov\v s -{\v S}{\v c}epin  for a finite dimensional $Y$ \cite{drrepsh})
Let $f: X \lo Y$ be a map of compacta such that $X$ is finite dimensional 
and $Y$ is full-valued  and let a CW-complex $K$ be such that 
$\edim (f^{-1}(y)\times Y) \leq K$ for every $y\in Y$. Then
$\edim X \leq K$.
\end{theorem}
{\bf Proof.} 
By Proposition \ref{2p0} and Theorem  \ref{dru} 
the theorem holds if either $Y$ or $f$  is $0$-dimensional.
If neither $Y$ nor $f$ is $0$-dimensional
then by Proposition \ref{2p2} we may assume that $K$ is simply connected.
Fix $i>1$ and assume  $H=H_i (K) \neq 0$.
Let $m=\dim_\z Y$, $n=\dim_H f$
and $k =$sup$\{ \dim_{H} (f^{-1}(y)\times Y): y\in Y\}$. By Theorem \ref{dr1},
$k \leq i$.
Since $Y$ is full-valued 
$\dim_{H} Y =m$ and $k = n +m$.
By (2') of Theorem \ref{t2}, Theorem \ref{dru} and Proposition \ref{2p0}
there exists
a $\sigma$-compact $A \subset X_*$ such that $\dim_{H} (A ) \leq  n$
and $f_* |_{X_* \sm A}$ is $0$-dimensional. Let $\sigma (H)$ be the Bockstein
basis for $H$ and let $L =\vee \{ M(G, n): G \in \sigma (H) \}$ where
$M(G,n)$ is a Moore space of type $(G,n)$.
Then $L$ is countable and by Bockstein theory and Theorem \ref{dr1}
 e-dim $A \leq L$.
By Theorem \ref{ol} there is a completion $A'\subset X_*$ of $A$ such that 
$\edim A'\leq  L$.
Then by Theorem \ref{dru}, $B'=X_* \sm A'$ is of $\dim =\dim_\z \leq  m$.

Then by Theorem \ref{dyd} $\edim X_* = \edim (A'\cup B')
\leq L*S^m= \Sigma^{m+1} L=
\vee \{ M(G,n+m+1): G \in \sigma(H) \}$. Hence by Theorem \ref{dr1}
$\dim_H X_* \leq n+m+1$. Then $\dim_H X \leq n+m =k$. Thus
$\dim_{H_i(K)} X \leq i$ and by Theorem \ref{dr1} the theorem follows.
\hfill $\Box$

Note that  one can avoid the use of  Olszewski's completion theorem
(Theorem \ref{ol}) and Bockstein theory in the proof of Theorem \ref{t4}
by using Theorem \ref{5t1} instead of Theorem \ref{dyd}.

\end{section}

\begin{section}{Krasinkiewicz maps}
Krasinkiewicz \cite{kr} introduced and studied maps having 
the following remarkable property. We will call 
these maps Krasinkiewicz maps.
\begin{definition}
(\cite{kr})
A map $f : X \lo Y$ is said to be a Krasinkiewicz  map if any continuum
in  $X$ either contains a component of a fiber of $f$ or is contained
in a component of a fiber of $f$.
\end{definition}

Krasinkiewicz \cite{kr} showed that for a compactum $X$ the Krasinkiewicz
maps form a dense subset in $C(X,I)$. We  improve this result
by showing that almost every map in $C(X,I)$ is a Krasinkiewiz map.

First we present a short construction of
Krasinkiewicz maps using a different  approach based on \cite{l2}.  
Namely, we will use the following proposition  actually proved in \cite{l2}.

\begin{proposition}
\label{lelek}
(\cite{l2}, Proposition 2.2, Proof B)
Let $Z \subset X$ be a closed $0$-dimensional subset of a compactum $X$.
Then for almost every map in $C(X,I)$, $Z$ is contained in the union 
of trivial components of the fibers of $f$.
\end{proposition}

We construct Krasinkiewicz maps as follows. Let $X$ be a compactum.  
 Take
a sequence of Cantor sets $C_i\subset I$ and
intervals $[a_i , b_i ] \subset I$, $i=1,2,...$ such that 
$C_i \subset (a_i, b_i)$, $[a_{i+1}, b_{i+1}] \cap (C_1\cup...\cup C_i)=\emptyset$
and each non-empty open
interval  in $I$ contains some $C_i$. 
Fix $\epsilon >0$ and let $f \in C(X,I)$.
Let $\psi_0 =f$ and assume that we have constructed $\psi_{i-1}: X \lo I$.
Let  $g_i : C_i \lo 2^X$ be a map onto
the  hyperspace $2^X$ of all non-empty closed subsets of $X$ equipped with
the Hausdorff topology. Define 
$A_i = \cup \{ g_i(c) \cap {\psi}_{i-1}^{-1}(c): c \in C_i \}$.
Then $A_i$ is closed in $X$. Let $q_i : X \lo X_i$ be the quotient map to
the space $X_i$ obtained from $X$ by identifying 
 each component of $A_i$ with
a singleton. Then $\psi_{i-1}$ factors through $\phi_i : X_i \lo I$ and 
$Z_i = q_i(A_i)$ is $0$-dimensional. 
By Proposition \ref{lelek} $\phi_i$ can be $\epsilon/2^{i+1}$-approximated
by a map $\phi'_i : X_i \lo I$ such that $Z_i$ is contained in the union
of trivial components of the fibers of $\phi'_i$. 
Let $\psi_i = \phi'_i \circ q_i : X \lo I$.
Then $\psi_i$ is $\epsilon/2^{i+1}$-close to $\psi_{i-1}$ and it easy to see that
we may assume that $\psi_i(x) = \psi_{i-1}(x)$ if 
$\psi_{i-1}(x) \in I\sm (a_i, b_i)$.
Then $\psi =\lim \psi_i$ is $\epsilon$-close to $f=\psi_0$ and one can easily
observe that $\psi$ is a Krasinkiewicz map. Moreover, $\psi$ satisfies the following
condition:

(*) for every continuum $F \subset X$ such that $\psi(F)$ is not a singleton
there is a subset $D \subset \psi(F)$ dense in $\psi(F)$ such that for every
$d\in D$ ,$\psi^{-1}(d) \cap F$ is the union of components of $\psi^{-1}(d)$.

\begin{proposition}
\label{kp2}
Let $X$ be a compactum.  The set of maps  in $C(X,I)$  satisfying 
the condition (*)  is  a dense $G_\delta$-subset of $C(X,I)$.
\end{proposition}
{\bf Proof.} Let $H_* =$the set of maps in $C(X,I)$ satisfying (*).
We have already shown that $H_*$ is dense in $C(X,I)$.
For $0\leq  p< q \leq 1$ 
denote by $H(p,q,n)$ the set of functions
$ f \in C(X,I)$ such that there is a continuum $F \subset X$ such that
$[p,q] \subset f(F)$ and for every $d \in [p, q]$ there is a component
of the fiber $f^{-1}(d)$ intersecting $F$ and containing a point $x$ such that
dist$(x, F) \geq 1/n$. It is easy to check that $H(p,q,n)$ is closed
in $C(X,I)$.  
Denote $H=\cup \{ H(p,q,n) : p,q $ are rationals, $0\leq p < q\leq 1, n=1,2,..\}$

Let us show that  
$H_* = C(X,I) \sm H$. It is clear that  $H_* \subset  C(X,I)\sm H$.
If $\psi$ does not satisfy (*) then there are a continuum $F \subset X$
and an interval $[a,b] \subset \psi(F)$ such that for every $d \in [a,b]$
there is a component of $\psi^{-1}(d)$ intersecting both $F$ and $X\sm F$.
Let $A_n = \{ d \in [a,b]: $ there are a component $C$ of $\psi^{-1}(d)$ 
 and a point $x \in C$ such that $C\cap F \neq \emptyset$ and
dist$(x, F) \geq 1/n \}$. Then $A_n$ is closed, 
$[a,b] = \cup A_n$  and hence there are $n$ and  a non-degenerate interval $[p,q]$ with
rational endpoints such that $[p,q] \subset A_n$. Then $\psi \in H(p,q,n)$
and we proved that $H_* =C(X,I)\sm H$.
Thus $H_*$ is $G_\delta$ in $C(X,I)$ and the proposition follows.
\hfill $\Box$

Clearly each map satisfying (*) is a Krasinkiewicz map and Proposition \ref{kp2}
implies:

\begin{theorem}
\label{kt1}
Let $X$ be a compactum.
Almost every map in $C(X,I)$ is a Krasinkiewicz map.
\end{theorem}

\begin{proposition}
\label{kp21}
Let $f : X \lo Y$ be  a Krasinkiewicz map of compacta $X$ and $Y$.
Denote $C(f)= $the union of all non-trivial components of the fibers of $f$.
Then there are compacta $A_1, A_2,... \subset X$ such that $C(f)=\cup A_i$
and each component of $A_i$ is contained in a fiber of $f$.
In particular $\edim C(f) \leq \edim f$.
\end{proposition}
{\bf Proof.} Let $C_n $ = the union of all non-trivial components
of $f$ of diam$ \geq 1/n$. Clearly $C_n$ is closed and $C(f)=\cup_{n=1,\infty} C_n$.
Cover $C_n$ by finitely many closed subsets $B^n_1,B^n_2,...\subset C_n$ of diam$ <1/n$.
Then a component of $B^n_j$ cannot  contain a component of a 
fiber of $f$ and hence must
be contained in a componet of a fiber of $f$. Enumerate $ \{ B^n_j \}$ 
as $\{A_1, A_2,...\}$
and we are done.
\hfill $\Box$

A continuum $X$ is said to be hereditarily indecomposable 
if for every pair
of intersecting subcontinua $A, B \subset X$, either $A \subset B$ or $B \subset A$.
A map $f$ is said to be a Bing map if each component  of each fiber of $f$ is
hereditarily indecomposable.

\begin{theorem}
\label{kt3}
(\cite{l1})
Let $X$ be a compactum.
Almost every map in $C(X,I)$ is a Bing map.
\end{theorem}
By a Bing-Krasinkiewicz map we mean a map which is both Bing and Krasinkiewicz.
By Theorems \ref{kt1} and \ref{kt3}
for a compactum $X$ almost every map in $C(X,I)$ is
Bing-Krasinkiewicz.
It is easy to verify:
\proposition
\label{kp3}
Let $f : X \lo Y$ be a Bing-Krasinkiewicz map of compacta.
 Then for every map $g : X \lo Z$ to a compactum $Z$
the map $(f,g) : X \lo Y \times Z$ is also Bing-Krasinkiewicz.

\end{section}

\begin{section}{Proofs of Theorems \ref{t1} and \ref{t2}}
We first present some facts used in the proofs.
\begin{proposition}
\label{pp1}
(\cite{st})
Let $Y$ be a compactum. Then there is a $\sigma$-compact
$0$-dimensional subset $A \subset Y\times I$
such that for the projection $p : Y\times I \lo Y$, $\dim p|_{(Y\times I) \sm A}=0$.
\end{proposition}
Such a subset $A$ is constructed as follows. Let $C_1,C_2,... \subset I$ be
a sequence of Cantor sets such that each non-empty open subset of $I$ contains some $C_i$.
Let a map $f_i : C_i \lo Y$ be onto. Define $A_i = \{(f_i(c), c):c\in C_i \}$.
Then $\dim A_i =0$ and  $A=\cup A_i$ satisfies Proposition \ref{pp1}.

\begin{proposition}
\label{pp2}
Let $X$ be a compactum, let $K$ be a CW-complex and let a map $f : X \lo I$
be of $\edim \leq K$. Then $\edim X \leq \Sigma K$.
\end{proposition}
{\bf Proof.} Let $g : F \lo \Sigma K$
be a map of a closed subset $F$ of $X$.
 Since $\edim f \leq K \leq \Sigma K$,
$g$ can be extended  over a small neighbourhood of  each fiber of $f$
and hence there is a partition $0=t_1 < t_2 <...< t_k=1$ of $I$ such that
$g$ extends to $g_i : F_i=F \cup f^{-1}([t_i, t_{i+1}]) \lo \Sigma K$
for every $i=1,k-1$. Since $\edim f \leq K$,
$\edim (f^{-1}(t) \times I) \leq \Sigma K$ for every $t \in I$ and hence
$g_{i} |_{f^{-1}(t_{i+1})}$ and $g_{i+1}|_{f^{-1}(t_{i+1})}$ are 
homotopic for every
$i=1,k-2$. It follows that there is a map $g' : X \lo \Sigma K$ which
is homotopic  on  $F_i$ to $g_i$ for each $i=1,k-1$
and hence $g$ can be extended over $X$.
\hfill $\Box$

\begin{theorem}
(\cite{dr1})
\label{pt1}
Let $K$ and $L$ be  countable CW-complexes and let for a compactum $X$,
$\edim X \leq K*L$. Then there is  a $\sigma$-compact $A \subset X$
such that $\edim A \leq K$ and $\edim (X\sm A) \leq L$.
\end{theorem}

\begin{proposition}
\label{pp3}
Let $f : X \lo Y$ be a light  map
of finite dimensional compacta $X$ and $Y$. Then
almost  every  map $g : X \lo I$ has the property
that each fiber of 
the map $(f,g) : X \lo Y\times I$ contains at most 
$\dim Y +1$ points;
\end{proposition}
{\bf Proof.}
Let $\delta >0$ and consider the set $\cal G$
of maps $g : X \lo I$ having the property
that each fiber of $(f,g)$ contains at most $k+1$ points
having pairwise distances $\geq \delta$.
Clearly $\cal G$ is an open subset
of the mapping space $C(X,I)$. We will show that $\cal G$ is also dense
proving the proposition.

Take a map   $g': X \lo I$ and
fix $\epsilon >0$. 
Let  $\phi_\alpha : Y \lo I$  be a finite set of maps
forming a partition of unity for $Y$
 such that $V_\alpha =$supp$\phi_\alpha =\phi^{-1}_\alpha (0, 1]$ 
 is a cover of order $k+1$.
Since   $\dim f=0$ we may 
 assume that each $f^{-1}(V_\alpha)$ splits into a finite family
of disjoint open sets $U_\alpha^1, U_\alpha^2,...$ such that 
diam $U_\alpha^i \leq \delta/3$ and diam$g'(U_\alpha^i)\leq \epsilon/3$
for every $i$. Let $\psi_\alpha^i : X \lo I$ be defined by
$\psi_\alpha^i (x)=\phi_\alpha (f (x))$ if $x \in U^i_\alpha$ and
$\psi_\alpha^i (x)=0$ if $x \in X \sm U^i_\alpha$.
Then $\psi_\alpha^i$ form a partition of unity for $X$. Choose real
numbers $0 < b_\alpha^i < 1$ such that dist$(b_\alpha^i , g'(U_\alpha^i)) < \epsilon/6$
 and $b_\alpha^i$  have the following
general position property: for every $l$
any collection of $l+1$ points in $R^{l}$ with the coordinates formed
from different elements $b_\alpha^i$ is not contained in a hyperplane 
of $R^{l}$.

Define $g : X \lo I $ by $g(x)=\Sigma_{\alpha, i} \psi_\alpha^i(x)b_\alpha^i$.
Then $|g(x)-g'(x)| \leq \epsilon$. We will show that $g \in \cal G$. Let 
 $a \in Y$,  $b \in I$
and let $x_1, x_2,... \in X$ be such that dist$(x_{j_1}, x_{j_2}) \geq \delta, j_1 \neq j_2$
and $(f,g)(x_j)=(a,b)$. Let  $V_{\alpha_1},...,V_{\alpha_l}, l\leq k+1$ be the
sets containing $a$.  Then $x_j \in U_{\alpha_t}^{i(j,t)}$ and $i(j_1,t) \neq i(j_2,t)$
if $j_1 \neq j_2$ for every $t=1,l$.  Now  $g(x_j)=
\Sigma_{t=1,l} \psi_{\alpha_t}^{i(j,t)}(x_j)b_{\alpha_t}^{i(j,t)}=
 \Sigma_{t=1,l} \phi_{\alpha_t}(f(x_j))b_{\alpha_t}^{i(j,t)}=
\Sigma_t \phi_{\alpha_t}(a)b_{\alpha_t}^{i(j,t)}=b$.
 Then for every $j$, $ (b_{\alpha_1}^{i(j,1)}, b_{\alpha_2}^{i(j,2)},...,b_{\alpha_l}^{i(j,l)})$
   is a solution of the linear equation         $\Sigma_t \phi_{\alpha_t}(a) s_t=b$
   for the variables $s_t $.  Since $ \Sigma_t \phi_{\alpha_t}(a)=1$   there is
   at least one non-zero coefficient in the equation and
 according to our choice of $b_\alpha^i$
there are  at most $l$  solutions of this equation
of the form    $ (b_{\alpha_1}^{i(j,1)}, b_{\alpha_2}^{i(j,2)},...,b_{\alpha_l}^{i(j,l)})$
  and hence
$j \leq l\leq k+1$.
Thus we have proved that $g \in \cal G$ and the proposition follows.
\hfill $\Box$
\\
{\bf Proof of Theorem \ref{t1}.}\\
{\bf (1)$\Rightarrow$(3)}
Take a finite-to-one map $\psi : C \lo Y$ from a Cantor set $C$
onto $Y$ and
let $Z= \{ (c,x) \in C\times X: \psi (c) =f(x)\}$
 be the pullback of  $\psi$ and $f$ with the projections 
$p_C : Z \lo C$ and $p_X : Z \lo X$.  Then $\edim p_C  =
\edim  f \leq \Sigma K$ and by Proposition \ref{2p0} $\edim Z \leq \Sigma K$.
By Theorem \ref{pt1} there are 
 $0$-dimensional compact subsets $A_1, A_2,...$ of $Z$  such that 
 $\edim(Z\sm A)\leq K$
where $A=\cup A_i$. Then $\edim(f^{-1}(y)\sm p_X( A)) \leq K$ 
for every $y\in Y$.
Since $\psi$ is finite-to-one we have that $f|_{p_{X}(A_i)}$ is 
$0$-dimensional for every $i$. Hence by Proposition \ref{pp3} 
 almost every map $ g : X \lo I$ is  
such that $(f,g)$ is at most $(\dim Y +1)$-to-$1$ on $p_X(A)$.
 Let $(y,t) \in Y\times I$.
Then $(f,g)^{-1}(y,t) \subset (f^{-1}(y) \sm p_X(A))\cup 
((f,g)^{-1}(y,t) \cap p_X(A))$
and hence $\edim(f,g) \leq K$.
 \hfill $\Box$
 \\
 {\bf (3)$\Rightarrow$(2)}
 By Theorems \ref{kt1} and \ref{kt3}  there is a Bing-Krasinkiewcz map 
 $g : X \lo I$ such that $\edim (f,g) \leq K$. By proposition \ref{pp1}
 there is a $\sigma$-compact $0$-dimensional set $S \subset Y \times I$
 such that $\dim ((\{ y \} \times I) \sm S)=0$ for every $y \in Y$.
 Then by Proposition \ref{2p0}, $A' =(f,g)^{-1}(S)$ is of $\edim \leq K$.
 By Proposition \ref{kp3} $(f,g)$ is Bing-Krasinkiewicz
 and by Proposition \ref{kp21}, $A''=C((f,g))$ is
 of $\edim \leq K$. 
 Then $A =A' \cup A''$ is $\sigma$-compact and we show
 that $f|_{X\sm A}$ is $0$-dimensional. Let $h_m : X  \lo Z$ and $h_l :
 Z \lo Y \times I$ be the monotone-light decomposition of $(f,g)=h_l \circ h_m$
 with $h_m$ monotone and $h_l$ light. Then for every $y \in Y$,
 $h_m |_{f^{-1}(y) \sm A''}$ is a homeomorphism and
 $h_m(f^{-1}(y) \sm A') \subset h_l^{-1}((\{ y \} \times I) \sm S)$.
 Since  $h_l^{-1}((\{ y \} \times I) \sm S)$ is $0$-dimensional and 
 $h_m(f^{-1}(y) \sm A) \subset h_l^{-1}((\{ y \} \times I) \sm S)$,
 $f^{-1}(y) \sm A$ is also  $0$-dimensional and 
 therefore $\dim f|_{X\sm A} =0$.
 \hfill $\Box$
 \\
 {\bf (3)$\Rightarrow$(3')} is obvious. \hfill $\Box$
 \\
 {\bf  (3')$\Rightarrow$(1)} by Proposition \ref{pp2}. \hfill $\Box$
 \\
 {\bf (2)$\Rightarrow$(1)}  by Theorem \ref{dyd}. \hfill $\Box$
 \\
 {\bf Proof of Theorem \ref{t2}.}\\
 {\bf  (0)$\Rightarrow$(3)} Let $g_*: X_* \lo I$ be a Bing  map.
 Then a fiber $(f_*,g_*)^{-1}(y,t)  = g_*^{-1} (t) \cap
 (f^{-1}(y) \times I ) , y\in Y, t\in I$
  contains no non-degenerate interval and hence the projection of
  $(f_*,g_*)^{-1}(y,t)$ onto $f^{-1}(y)$ is a  $0$-dimensional map.
  Thus by Theorem \ref{dru} $\edim (f_*,g_*)^{-1}(y,t) \leq \edim f^{-1}(y)$
  and hence $\edim (f_*, g_*) \leq K$. Note that by Theorem \ref{kt3}
  almost every map is a Bing map and we are done. \hfill $\Box$
  \\
  {\bf (3)$\Rightarrow$(2')} Replace $X,f$  and $g$ by $X_*, f_*$ and $g_*$
  respectively
  and use the construction and notation
   of the proof (3)$\Rightarrow$(2) of Theorem \ref{t1}.
   We only need to check that $A$ can be decomposed
   into a countable family of compacta whose components admit
    $0$-dimensional maps into the fibers of $f$. 
    Each  continuum contained in  $A'$ is contained in a fiber of $(f_*,g_*)$.
     By Proposition \ref{kp21}
    $A''$ can be decomposed into a countable
   family of compacta whose components are contained in the fibers
of $(f_*,g_*)$. Each continuum in a fiber of $(f_*, g_*)$ is hereditarily
indecomposible and hence its projection to the corresponding fiber
of $f$ is a $0$-dimensional map. \hfill $\Box$
\\
{\bf (2')$\Rightarrow$(2)} is obvious. \hfill $\Box$
\\
{\bf (2)$\Rightarrow$(1)} follows from Theorem \ref{dyd}. \hfill $\Box$
\\

\end{section}

\begin{section}{Remarks}
{\bf Remark 1.}
The properties (1),(2) and (3) of Theorem \ref{pt} are equialent
to:

(4) (cf. \cite{p1})
 $almost$  $every$  $map$ $g : X \lo I^{n+1}$ $has$ $the$ $property$
$that$ $each$ $fiber$ $of$ 
$the$ $map$ $(f,g) : X \lo Y\times I^{n+1}$ $contains$ $at$ $most$ 
$\dim Y+n+1$ $points$.

Indeed, by Proposition \ref{pp3} and Theorem \ref{pt}
this property holds if $\dim f \leq n$. Assume that there is $y \in Y$
such that $\dim f^{-1}(y) \geq n+1$. Take disjoint closed  
subsets $F_1, F_2,...,F_k \subset f^{-1}(y)$,
 $k=\dim Y +n +2$ such that $\dim F_i \geq n+1$
and let $g : X \lo I^{n+1}$ be such that $g|_{F_i} $ is essential
for each $i$ (that is $g|_{{g}^{-1}(\partial I^{n+1} )\cap F_i}$ cannot
be extended over $F_i$ as a map to $\partial I^{n+1}$). Then
any sufficiently close approximation of $g$ must be at least
$(\dim Y +n +2)$-to-$1$ on $f^{-1}(y)$
 and this contradiction shows that (4) implies
that $\dim f \leq n$.
\\\\
{\bf Remark 2.} 
The proof of Theorem \ref{dyd} given
in \cite{dyd} applies to prove the following stronger result.
\begin{theorem}
\label{5t1}
Let $K$ and $L$ be CW-complexes and 
let a separable metric space $X$ be decomposed into subsets
$X =A \cup B$ such that $\edim A \leq K$ and for every subset $F$
closed in $X$ and contained in $B$, $\edim F \leq L$. Then
$\edim X \leq K*L$.
\end{theorem}
Theorem \ref{5t1} allows one to avoid the use of Olszewski's
completion theorem in Theorem \ref{2t1} and  to extend it to the following
conclusion omitting the requirement that $K$ is countable.
\begin{theorem}
\label{5t2}
Let $K$ and $L$  be  CW-complexes such that for a map of compacta
$f : X \lo Y$, $\edim f \leq K$ and $\edim Y \leq L$. Then
$\edim X_* \leq K*L$ where $X_*=X\times I$.
\end{theorem}

E. {\v S}{\v c}epin 
conjectured that $\edim (X_1 \cup X_2) \leq \edim (X_1 * X_2)$.
This would significantly improve Theorem \ref{dyd}. Let us state without
a proof
the following result related to {\v S}{\v c}epin's conjecture.
\begin{theorem}
\label{5t3}
Let $X=X_1 \cup X_2$ be a decomposition of a compactum $X$.
Then $\edim X \leq \edim (\beta X_1 * \beta X_2)$, where $\beta X$ is the Stone-Cech
compactification of $X$.
\end{theorem}

Theorem \ref{5t3} implies, for example, that 
dim$X \leq \dim (\beta X_1 \times \beta X_2) +1$. 
Note that Theorem \ref{5t3} does not seem to be useful for infinite 
dimensional spaces, see \cite{l4}.
\\\\
{\bf Remark 3.} 
It seems to be of interest to know some dimensional
properties of $B = X_* \sm A$ in Theorem \ref{t2}, especially, when
$K$ is uncountable and  it is unknown if $A$ has a completion of the 
same extensional dimension. In general, Theorem \ref{dru} does not hold
for spaces  which are not $\sigma$-compact
and hence the $0$-dimensionality of $f_* |_B$
does not give much information about $B$.  In view of Theorem \ref{unif}
one is tempted to replace the $0$-dimensionality  by the uniform
 $0$-dimensionality.

Unfortunately 
Theorems \ref{pt}, \ref{t1} and \ref{t2} do not hold if the $0$-dimensionality
is replaced by the uniform $0$-dimensionality. Indeed, let $ p : I^n \lo I$, $n\geq 2$ be
the projection $p(x_1,...,x_n)=x_n$. Then for every $\sigma$-compact $(n-2)$-dimensional
subset $A$ of $I^n$ and every interval $(a,b)$ in $I$, $p^{-1}((a,b))\sm A$ is connected.

However one can observe that  the following result is contained
in the proof of Theorem \ref{t2}. 

\begin{theorem}
Let $f : X \lo Y$ be a map of compacta
 and let $f_* : X_* =X \times I \lo Y$ be
defined by $f_*(x,t)=f(x)$.
Then for almost every map $g_* : X_* =X\times I \lo I$
there is a $\sigma$-compact $A\subset X_*$ such that $A$ splits
into a countable family of compacta whose components admit
 $0$-dimensional maps into the fibers of $f$ and for $B = X_* \sm A$,
 $ f_* |_{B} $
is $0$-dimensional  and 
 $(f_*,g_*)|_{B}$ is uniformly 
$0$-dimensional.
\end{theorem}

In particular $\edim A \leq \edim f$ and $\edim  B \leq \edim \Sigma Y$.
\\\\
{\bf Remark 4.}
 Let us repeat that it is unknown if Theorems \ref{pt}
and \ref{t1} hold without any dimensional restriction on $Y$. We end this note 
with posing the following problem.
\begin{problem}
(cf. Theorem \ref{pt1}) Let $f : X \lo Y$ be a map of compacta
with $Y$ finite dimensional and let $K$ and $L$ be countable
CW-complexes such that  $\edim f \leq K*L$. Does there exist
a  $\sigma$-compact $A \subset X$ such that $\edim A \leq K$ and
$\edim f|_{X \sm A} \leq  L$ ?
\end{problem}

\end{section}

Department of Mathematics\\
Ben Gurion University of the Negev\\
P.O.B. 653\\
Be'er Sheva 84105, ISRAEL  \\
e-mail: mlevine@math.bgu.ac.il\\\\
Department of Mathematics\\
Texas Tech University\\
Lubbock, TX 79409-1042, USA\\
e-mail:   wlewis@math.ttu.edu

\end{document}